\documentclass{article}                                      
\usepackage{amsmath,amssymb,amscd,amsfonts} 
\usepackage{graphicx} 
\newtheorem{theorem}{Theorem}[section] 
\newtheorem{lemma}{Lemma}[section]
\newtheorem{observation}{Observation}[section]

 \newcommand{\proof}{\noindent{\it Proof:~}}
\newcommand{\KONEC}{{\hfill$\square$}}
\newcommand{\eps}{\varepsilon}
 
\makeatletter
\newcommand*{\bdiv}{%
  \nonscript\mskip-\medmuskip\mkern5mu%
  \mathbin{\operator@font div}\penalty900\mkern5mu%
  \nonscript\mskip-\medmuskip
}
\makeatother

\begin{document}

\title{\bf  Improved approximative multicoloring of hexagonal graphs }

%\small
{
\author{ 
Janez \v{Z}erovnik\thanks{Supported in part by ARRS, the research agency of Slovenia.}\\
FME, University of Ljubljana,\\
%Faculty of Mechanical Engineering \\
A\v sker\v ceva 6, % \\
SI-1000 Ljubljana, Slovenia\\
{\tt janez.zerovnik@fs.uni-lj.si	}
}
\date{\today}

\maketitle

\begin{abstract}  
In 1999,   McDiarmid and   Reed    conjectued that 
  the approximation ratio $9/8$ of multichromatic number to weighted clique number asymptotically is  the best possible for general weighted  hexagonal graphs. 
 We prove that 
 there is a proper multicoloring of $G$ that uses at most   $15 \lfloor\frac{\omega(G)}{12}\rfloor+ 18$     colors 
  improving the best previously  known asimptotic ratio from 4/3  to 5/4.
\end{abstract}

\section{Introduction}

A fundamental problem concerning cellular networks is to assign sets of frequencies (colors) to transmitters (vertices) in order to avoid unacceptable interferences~\cite{Hale}. 
The number of frequencies' demands from a transmitter may vary between transmitters. 
In a usual cellular model, transmitters are centers of hexagonal cells and the corresponding adjacency graph is a subgraph of the infinite triangular lattice. 
An integer $d(v)$ is assigned to each vertex of the triangular lattice and will be called the \emph{weight} (or \emph{demand}) of the vertex $v$. 
The vertex weighted graph induced on the subset of vertices of the triangular lattice is called a \emph{hexagonal graph}, and is denoted by $G=(V,E,d)$. 
%Hexagonal graphs arise naturally in studies of cellular networks. 
A \emph{proper multicoloring} of $G$ is a mapping $f$ from $V(G)$ to subsets of integers such that $\left\vert f(v)\right\vert \geq d(v)$ for any vertex 
$v\in V(G)$ and $f(v)\cap f(u)=\emptyset$ for any pair of adjacent vertices $u$ and $v$ in the graph $G$. 
The minimal cardinality of a proper multicoloring of $G$   is called the \emph{multichromatic number}, 
denoted   $\chi (G)$   (for simplicity,   as it is the only chromatic number studied here). %%%%%%%%%%%%%%%%%% $\chi_{m}(G)$,
Another invariant of interest in this context is the \emph{(weighted) clique number}, $\omega(G)$, defined as follows: 
The weight of a clique of $G$ is the sum of weights  on its vertices and $\omega(G)$ is the maximal clique weight on $G$. Clearly, $\chi(G) \geq \omega(G)$ and a
simple example shows that the upper bound is 
 $\chi(G)\leq(4/3)\omega(G)$ (see Figure  \ref{exampleW3_X4}).
However, note that in the example,   $\omega(G) = 3$ and $\chi(G)=4$
and we are not aware of any example with  $\omega(G) > 3$ and $\chi(G)=4 \omega(G) /3 $.
Interesting enough, 
the best known approximation algorithms for hexagonal graphs give solutions within  $\chi(G)\leq(4/3)\omega(G)+O(1)$ \cite{pierwsza,cykle,a43}.
For a later reference note that the constant is very small, 
more precisely, we have $\chi(G)\leq \frac{4\omega(G)+1 }{3} \leq  4\lceil\frac{\omega(G)}{3}\rceil$.

Better bounds have been obtained for triangle-free hexagonal graphs. 
In \cite{havetJZ} a distributed algorithm with competitive ratio 5/4 is given. 
Existence of 2-local distributed algorithm with competitive ratio 5/4 is  provided  in   \cite{a54}, and
  an inductive proof for ratio 7/6 is reported in \cite{Havet}. A 2-local 7/6-competitive algorithm for a sub-class of triangle-free hexagonal graphs is given in \cite{a76subclass}. 
As the shortest odd cycle that can be realized as an induced subgraph of triangular lattice has length 9, the 
best possible approximation      is  $\chi(G)\leq (9/8)\omega(G) + O(1)$.
A  challenging research  question is to settle the  conjecture due to McDiarmid and Reed:

\medskip
\noindent{\bf Conjecture}   \cite{pierwsza}  There is an absolute constant $ c $ such that for every weighted hexagonal graph $ G $,
$$\chi(G) \leq \frac{9}{8}\omega(G) + c $$ 

The conjecture remains unsolved even in the  relaxed version for triangle-free hexagonal graphs. 

An interesting  special case of a proper multicoloring is when $d$ is a constant function. 
In this case we  define a $n$-$[k]$-coloring to be an assignment of sets of $k$ colors from a set of 
  $n$ colors to each vertex.
  For example a $9$-$[4]$-coloring is an assignment of four colors from the set $\{1, ... ,9\}$ to each vertex. 
In \cite{havetJZ},  a  $5$-$[2]$-coloring    algorithm is used to provide a general 5/4-competitive algorithm.
An elegant idea that implies the existence of a $14$-$[6]$-coloring   
was presented in~\cite{tfree}.
The existence of a $7$-$[3]$-coloring  follows from result of \cite{Havet}.
A shorter proof  based on the 4-colour theorem  is provided in \cite{ost}.
In \cite{a73lin} a linear time algorithm for $7$-$[3]$-coloring  of an arbitrary triangle-free hexagonal graph $G$ is given. 
The only  result regarding $9$-$[4]$-coloring, however restricted to a  subclass of hexagonal graphs, is the algorithm for
 $9$-$[4]$-coloring of triangle-free hexagonal graphs without neighboring corners \cite{Utilitas}.
The $n$-$[k]$-colorings can be used to obtain $\frac{n}{2k}$ approximations for $\chi(G)$ 
due to   the following lemma.
(Folklore, for example the proof for $k=4$ appears in \cite{Utilitas}.)

\begin{lemma}\label{eqivalention}
For each triangle-free graph $G$ the following statements are equivalent:
\begin{itemize}
 \item[(i)] $G$ is $(2k+1)$-$[k]$-colorable, 
 \item[(ii)] $G$ is $\left\lceil \frac{2k+1}{2k}\omega(G) \right\rceil$ multicolorable.
\end{itemize}
\end{lemma}
 
On the other hand,  triangle-free hexagonal graphs with large odd girth surely allow better approximation, 
since hexagonal graphs are planar. 
However, we are not aware of any work on  hexagonal graphs that would  improve 
the well known  results that hold for arbitrary planar graphs  \cite{Klostermeyer1,Klostermeyer2}. 
Another  avenue that may be of interest is   to study  generalization of the (planar) hexagonal graphs
to  3D, where the triangular grid is replaced by the "cannonball" grid 
 \cite{clanek3D}.

Here we consider the general hexagonal graphs, and concentrate on improvement of  the approximation 
bound. The basic idea is to partition the original hexagonal graph into three triangle-free hexagonal graphs 
and apply the algorithms for triangle-free graphs. Roughly speaking, assuming that the weights can be evenly 
divided among  the three subgraphs, then any approximation bound for the triangle-free hexagonal graphs 
implies the same approximation bound for the general case. 
The main new result  presented here is Theorem \ref{theorem54}
implying  that  there   is    an algorithm that  
colors any weighted  hexagonal graph   with $\chi(G) \leq \frac{5}{4}\omega(G) + O(1) $  colors.
This significantly  improves the previously best known bound  $\chi(G) \leq \frac{4}{3}\omega(G) + O(1) $ \cite{pierwsza,cykle,a43}.
The mail idea   might  be used to improve the bound to get closer to the conjectured 
 $\chi(G) \leq \frac{9}{8}\omega(G) + O(1) $. %\cite{pierwsza}.

The rest of the  paper is organized as follows.
In the next section, we formally define some basic terminology. 
The general method and its application using a 5/4 approximation algorithm providing the new bound is given in Section 3.
The last section provides a short conclusion and discusses ideas for  future work. 

\begin{figure}[htp] 
  \begin{center}
%   \leavevmode
  \includegraphics[height=3.5cm]{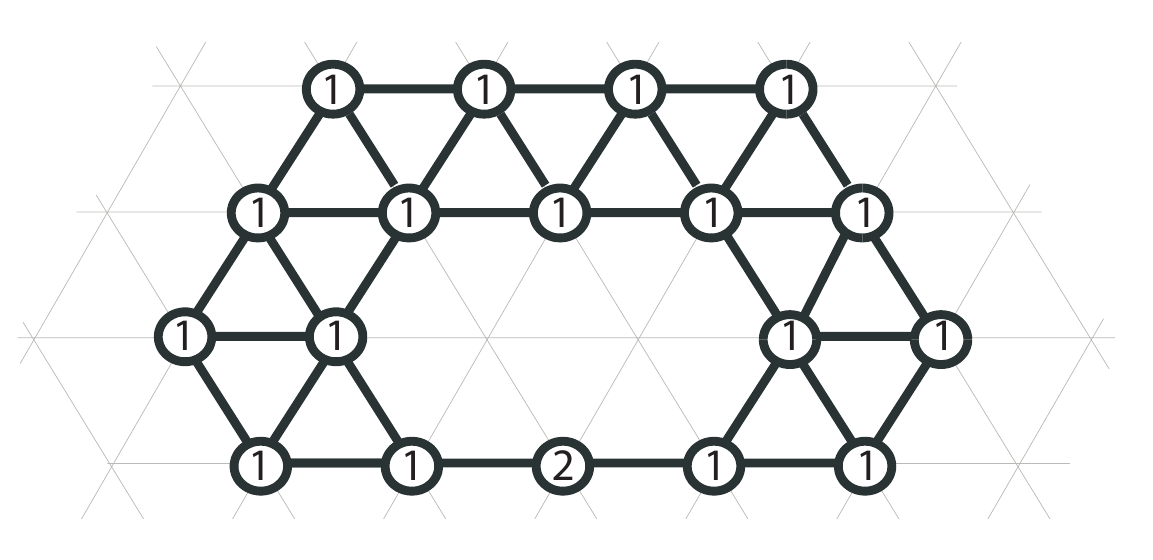} 
 \end{center}
\vskip -8mm
 \caption{An example of a  hexagonal graph.} 
 \label{exampleW3_X4}
\end{figure}

\section{Basic definitions and useful facts}\label{preliminary}

A vertex weighted graph is given by a triple $G(E,V,d)$, where $V$ is
the set of vertices, $E$ is the set of edges and $d: V \to \mathbb{N}$ is a weight function assigning    integer weights  to vertices of $G$.

The  vertices of the triangular lattice $T$ can be described as follows: 
the position of each vertex is an integer linear combination $x\vec{p}+y\vec{q}$ of two vectors $\vec{p}=(1,0)$  and $\vec{q}=(\frac{1}{2}, \frac{\sqrt{3}}{2})$. Thus vertices of the triangular lattice may be identified with pairs $(x,y)$ of integers. Given the vertex $v$ we will refer to its coordinates as $x(v)$ and $y(v)$.
Two vertices are adjacent when the Euclidean distance between them is one. 
Therefore each vertex $(x,y)$ has six neighbors: $(x-1, y), (x-1,y+1), (x,y+1),$ $(x+1,y), (x+1,y-1), (x,y-1)$. 
For convenience  we refer to the neighbors as: {\em left}, {\em up-left}, {\em up-right}, {\em right}, {\em down-right} and {\em down-left} neghbor.
Assume that we are given a weight function $d: V \to \{0,1,2, \ldots \}$ on vertices of triangular lattice.
We define a~{\em weighted hexagonal graph} $G=(V,E,d)$ as an induced subgraph on vertices of positive demand on the triangular lattice (see Figure \ref{hexagonal}).
Sometimes we want to work with (unweighted) hexagonal graphs $G=(V,E)$  that can be defined as induced graphs on the subset of vertices of the triangular lattice. In such a graph we can say that the weight function in each vertex in the graph is equal to 1, and for vertices out of the graph it is equal to 0.  

\begin{figure}[h]
 \centering
  \leavevmode
  \includegraphics[height=4.5cm]{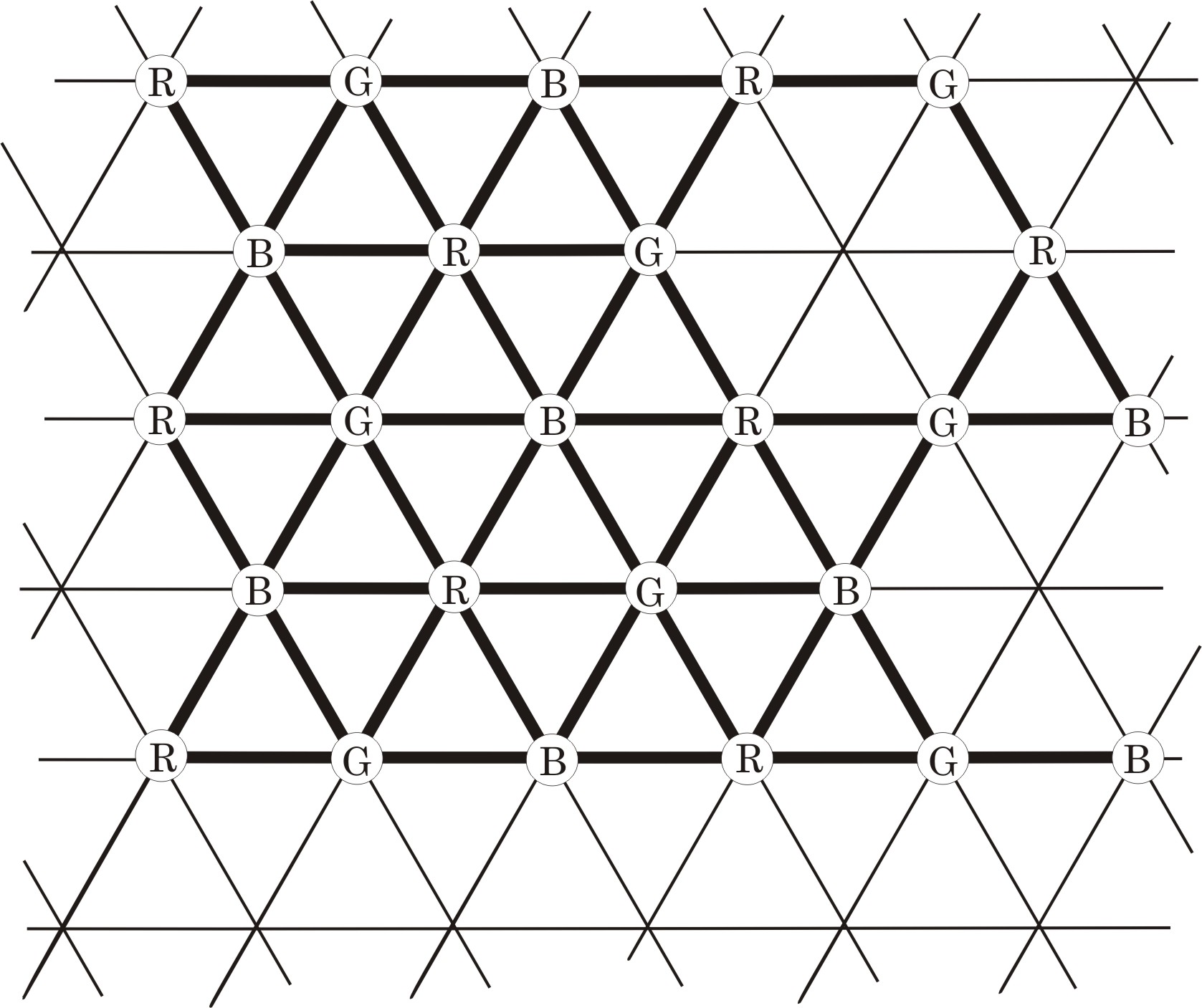}  
\hfil
  \includegraphics[height=4cm]{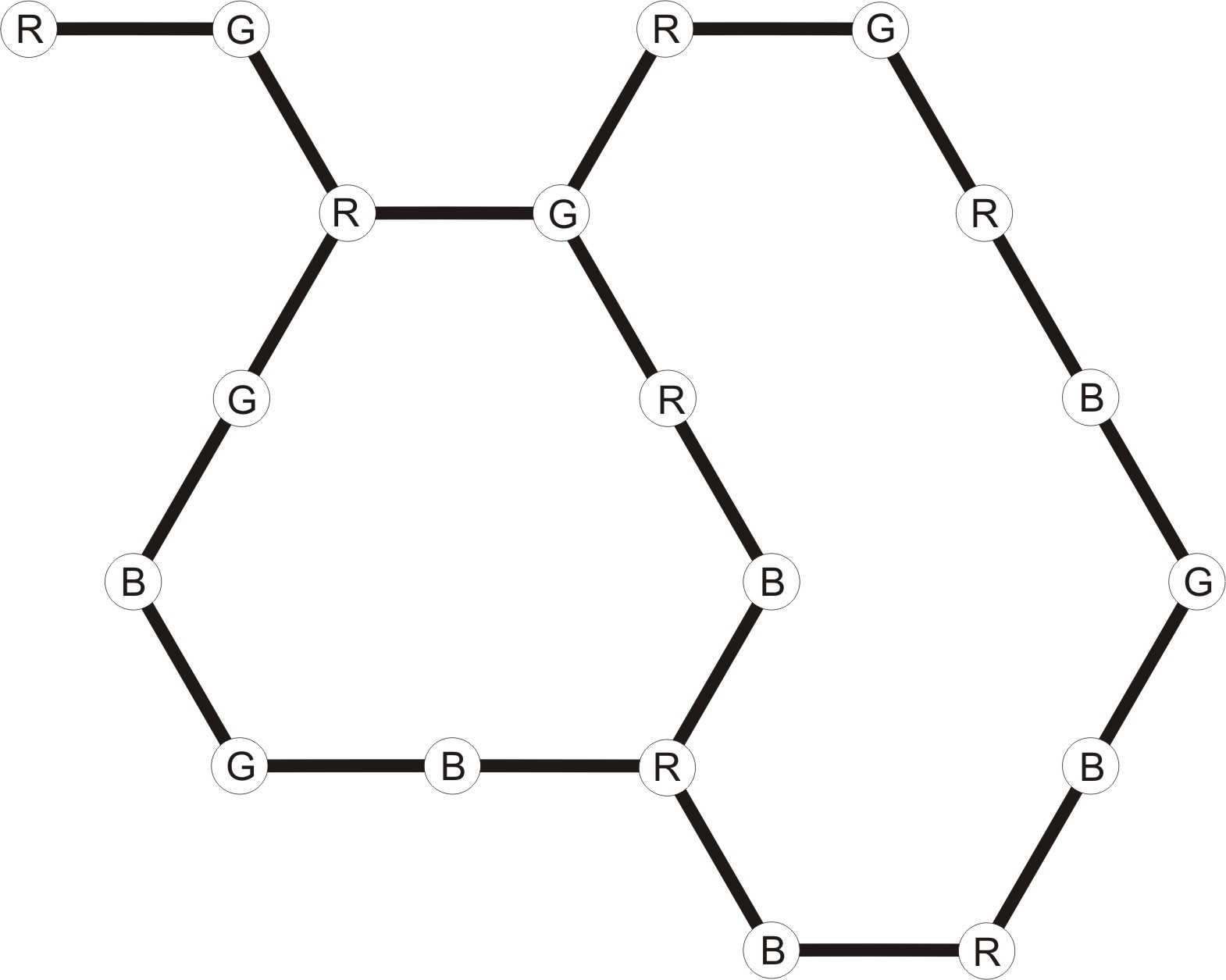}
 \caption{An example of a~hexagonal graph and a triangle-free hexagonal graph  (with base colorings).} 
 \label{hexagonal}
  \label{tfreehexagonal}
\end{figure}

There exists an obvious 3-coloring of the infinite triangular lattice which gives the partition of the vertex  set of any hexagonal graph into three independent sets. 
Let us denote the color of vertex $v$ in this 3-coloring by $bc(v)$ and call it a~{\em base color} (for simplicity we will use {\em red}, {\em green} and {\em blue} as base colors and their arrangement is given in Figure \ref{hexagonal}), 
i.e. $bc(v) \in \{ red, blue, green\}$.

We call a~{\em triangle-free hexagonal graph} an induced subgraph of the triangular lattice without  a  3-clique (see Figure~\ref{tfreehexagonal}).

%%%%%%%%%%%%%%%%%%%%%%%%%%%%%%%%%%%%%%%%%%%%%%%%%%%%%%%%%%%%%%%%%%%%%%%%
%
\section{The general  method and the 5/4 approximation} 
%
%%%%%%%%%%%%%%%%%%%%%%%%%%%%%%%%%%%%%%%%%%%%%%%%%%%%%%%%%%%%%%%%%%%%%%%%

The main idea is to partition  the weighted  hexagonal graph into three triangle-free hexagonal (sub)graphs  such that 
the weights  are  divided among the partition so that the 
clique number of each subgraph is about one third of the original clique number. 
More formally, given a weighted hexagonal graph $G(V,E,d)$, 
the idea is    to construct triangle-free subgraphs
$G_1(V_1,E_1,d_1)$, 
$G_2(V_2,E_2,d_2)$,  and 
$G_3(V_3,E_3,d_3)$, 
such that $G_i$ is induced on $V_i$ (for $i\in \{1,2,3\}$), 
and $d_1(v) + d_2(v) + d_3(v) = d(v)$ for any $v\in V$.
Among the decompositions we prefer the decompositions for which 
$\omega(G_1) + \omega(G_2) + \omega(G_3)  = \omega(G)$, 
or, at least 
$\omega(G_1) + \omega(G_2) + \omega(G_3)  = \omega(G) +C$ for some absolute constant $C$.
Having this, the  second step would be  to  multicolor  with  an approximation algorithm each  of the subgraphs that were constructed, 
therefore obviously we would  obtain the same approximation ratio.

At present, we are not able to provide such a decomposition in general case. 
However, below we show how to  decompose a special class of weighted hexagonal graphs
which in turn enables an iterative method that achieves 5/4 approximation for general hexagonal graphs.
%
%In the next sections we wil elaborate this idea to show that any algorithm that gives  5/4-approximation  of triangle-free hexagonal graphs  can be used to obtain a 
%5/4-approximative multicoloring of arbitrary hexagonal graph. 

%%%%%%%%%%%%%%%%%%%%%%%%%%%%%%%%%%%%%%%%%%%%%%%%%%%%%%%%%%%%%%%%%%%%%%%%
%
In the rest  of  this section, we will first provide a decomposition of 
a special class of  weighted hexagonal graph into three triangle-free weighted hexagonal graphs.
For these graphs, clique number is 12 by definition, and we show how to multicolor them with at most 15 colors.
Finally, by iterated application of the decomposition that reduces the clique number we prove the main result,  Theorem \ref{theorem54}.

%%%%%%%%%%%%%%%%%%%%%%%%%%%%%%%%%%%%%%%%%%%%%%%%%%%%%%%%%%%%%%%%%%%%%%%%
%
\subsection{ A special class of  weighted hexagonal graphs}   \label{partition}
%Weights of vertices in  subgraphs $G_{RB}$, $G_{BG}$,  and   $G_{GR}$  - special case } 
%
%%%%%%%%%%%%%%%%%%%%%%%%%%%%%%%%%%%%%%%%%%%%%%%%%%%%%%%%%%%%%%%%%%%%%%%%

%%%%%%%%%% OZNACCI TA GRAF z npr.  

\newcommand{\GG}{{\tilde G}}
\newcommand{\VV}{{\tilde V}}
\newcommand{\EE}{{\tilde E}}
\newcommand{\dd}{{\tilde d}}

For a  weighted hexagonal graph  $\GG = (\VV,\EE,\dd)$  considered in this section we assume   that 
the vertices of $\GG$   have weight $\dd(v)=4$ when $v$ is on a triangle.
 A vertex that is not on any triangle has   weight $\dd(v)=4, 6,$ or $8$   so that  the weight on any clique is at most 12.
(Hence, any neighbor of a vertex of weight 8 must have weight 4,    while a neighbor of a vertex of weight 6 may have weight 6 or 4.) 
Finally, if $v$   is an   isolated vertex, then $\dd(v) = 12$.
Say a vertex $v$ is  \emph{ light}   if $\dd(v) < 6$.

%
%Clearly, the clique number $\omega(\GG) = 12$, and furthermore, $\omega$  is a maximal 
%assignment   of weights   (i.e. one of the maximal assignments)  on  $(\VV,\EE)$ such that $\omega(\GG) = 12$. 

Construct the subgraphs     $G_{RB}$, $G_{BG}$,  and  $G_{GR}$   of $\GG$ as follows.
Recall the basic 3-coloring. 
The subgraph $G_{RB}$ is   induced on the vertex set $V_{RB}$ given by the following two rules:
(1) If a  vertex $v \in V(G)$   lies in a triangle, and is of base color {$R$}  or {$B$}, then $v\in V_{RB}$.
(2) If a   vertex $v$    has weight at least 6,  then $v\in V_{RB}$.
Or, in other words: $V_{RB}$ is $V(\GG)$ without  light  vertices of base color green    $G$. 
% (We say a vertex is light if it has demand less than $\omega(G)/3$.)
%
Clearly,  $G_{RB}$ is triangle-free.

Triangle-free subraphs $G_{BG}$ and   $G_{GR}$  are defined analogously.
(They are subgraphs of $\GG$ induced on the vertex set  of $\GG$ without  light vertices of  the corresponding base color.)

\begin{observation}\label{observation} 
Assume the  subgraphs    $G_{RB}$, $G_{BG}$, and   $G_{GR}$  of  $\GG$ are constructed as explained above.
Then the following is true:
(1) Each vertex  of $G$ that lies on a triangle  appears  in exactly two subgraphs.
(2) Each vertex that is not light  and hence  does not lie on a triangle is in all 3 subgraphs.
\end{observation}

\noindent 
Define the weights for  each  $G_\star \in \{G_{RB}, G_{BG},  G_{GR}\}$   as follows: 
\begin{itemize}
\item
Any isolated vertex in $G_\star$   gets weight 4.
\item 
A vertex of $G_\star$ that lies on a triangle in $\GG$  gets weight 2.  
\item 
A vertex $v$ that is not on a triangle in $\GG$  is assigned weight
     \begin{itemize}
          \item   $\dd(v) = 2$ if it has more than one  neighbor in $G_\star$,
          \item   $\dd(v) = 3$ if it has  exactly one   neighbor in $G_\star$, and 
          \item   $\dd(v) = 4 $ if it  is isolated  in $G_\star$.
     \end{itemize} 
\end{itemize}

It is easy to see that for any vertex $v\in \VV$, the sum of weights assigned to $v$ in   $  G_{RB}, G_{BG},$ and $G_{GR}$  is at least $\dd(v)$. 
 
Some cases are trivial, namely:
If   $v$ is on a triangle in $\GG$, then it is assigned demand 2 in each of the subgraphs, and  altogether $2+2+2  = 6 =\dd(v)$ as needed.
If  $v$ is isolated in $\GG$, then it is isolated in each of the subgraphs, and hence $4+4+4  = 12 =\dd(v)$.

If   $v$ is not on a triangle and is not isolated in $\GG$, then it may have weight 4, 6, or 8.
A light vertex ($\dd(v)=4$)  appears in two subgraphs, and  it gets 2+2=4 colors.
A vertex of weight 6 appears in all three subgraphs,  and  therefore  receives  2+2+2=6,  which  is sufficient.
The case when  $\dd(v)=8$ is settled by the following lemma
showing that in this case the demand 8  can always be partitioned, either as 4+2+2 or as 3+3+2.

\begin{lemma}\label{lemmaHeavyVertices} 
If a vertex of $\GG$ has weight 8,  then it is either isolated in one of the subgraphs, or it is a   vertex  of degree one in at least two of the subgraphs.
\end{lemma}

%\begin{proof}  
\proof
Recall that weight 8 is assigned to  vertices that are not on any triangle, and have only neighbors of weight 4.
Consider a vertex $v$ with weight 8, and  wlog assume it is of base color red.
As $v$ is not on a triangle, it may have degree at most 3.

First,  assume the vertex of weight 8 has exactly two neighbors of different base colors
(see  Figure   \ref{SpecialHeavy} (a)).
By definition, the green neighbor is removed in the graph $G_{RB}$ and hence $v$ has only one neighbor in the graph $G_{RB}$.
Similarly, the blue neighbor is removed in   the graph $G_{GR}$ and hence $v$ has only one neighbor in the graph $G_{GR}$.
Therefore $v$ is a vertex of degree one  in two of the subgraphs as claimed.

The second case is when all neighbors of $v$ are of the same base color. 
There can be one 
(Figure   \ref{SpecialHeavy} (b)), two, or three neighbors   (Figure   \ref{SpecialHeavy} (c)).
Wlog say the neighbors are of base color blue.
Then  $v$  is isolated  in the graph $G_{GR}$,  because  the blue neighbors are removed in   the graph $G_{GR}$.
\KONEC

\begin{figure}[htp] 
  \begin{center}
%   \leavevmode
  \includegraphics[height=2.5cm]{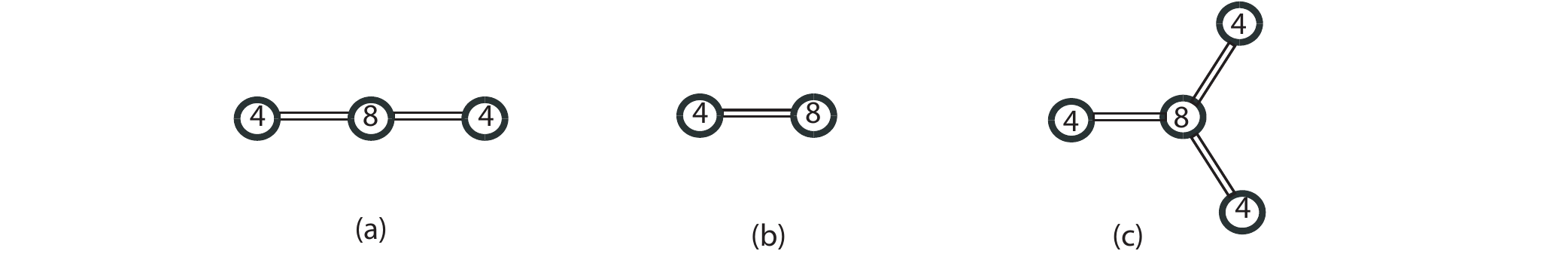}
 \end{center}
 \caption{ A configuration in which the heavy vertex  ($w=8$)  has  two neighbors of different base colors  (a) 
and  two configurations in which the heavy vertex has neighbors of the same base color (b) and (c). } 
 \label{SpecialHeavy}
\end{figure}

\begin{lemma}\label{lemma12}
There is a proper multicoloring of $\GG$ with 15 colors. 
\end{lemma}

%\begin{proof}
\proof
First, use   an  algorithm  for 5-[2] multicoloring   \cite{havetJZ} to   assign two colors to each vertex of demand 2 of   the graphs   $G_{RB}$, $G_{BG}$,  and   $G_{GR}$. 
(Use three distincts sets of  five colors.)
This  assigns  four  colors to each vertex of $\GG$ that lies on a triangle   and    up to six colors to vertices that appear in all three subgraphs
(in particular, it provides exactly six colors to vertices that have demand 2 in all three subgraphs). 

Next, observe that vertices of demand 3 in  $G_\star$  have at most one neighbor, so there are 3 free colors to assign.

Finally, note that   the vertices of demand 4 are isolated    in  $G_\star$,   hence they  can trivially  be assigned any set of   4   colors (out of  available 5 colors).

Hence, each of the graphs   $G_{RB}$, $G_{BG}$,  and   $G_{GR}$   is  properly colored using 5 colors, 
and consequently $\GG$  we have  a proper coloring of   $\GG$ with 15 colors.
\KONEC

\noindent
{\bf   Example.}
 Consider the graph on Figure    \ref{Example_decomposition}(a), 
and three subgraphs on \ref{Example_decomposition}(b),  (c) and (d).
If, for example, the  base color $bc(v)$ is red, then the subgraphs (b), (c) and (d) are   $G_{GR}$,   $G_{BG}$,  and   $G_{RB}$, respectively.
The heavy vertex $v$ has degree one  in two of the subgraphs (b) and (d), and hence it can be assigned three colors in each of these two subgraphs.
Thus we have a proper multicoloring using 15 colors. 

\begin{figure}[thhhb] 
  \begin{center}
  %  \leavevmode
  \includegraphics[height=9.5cm]{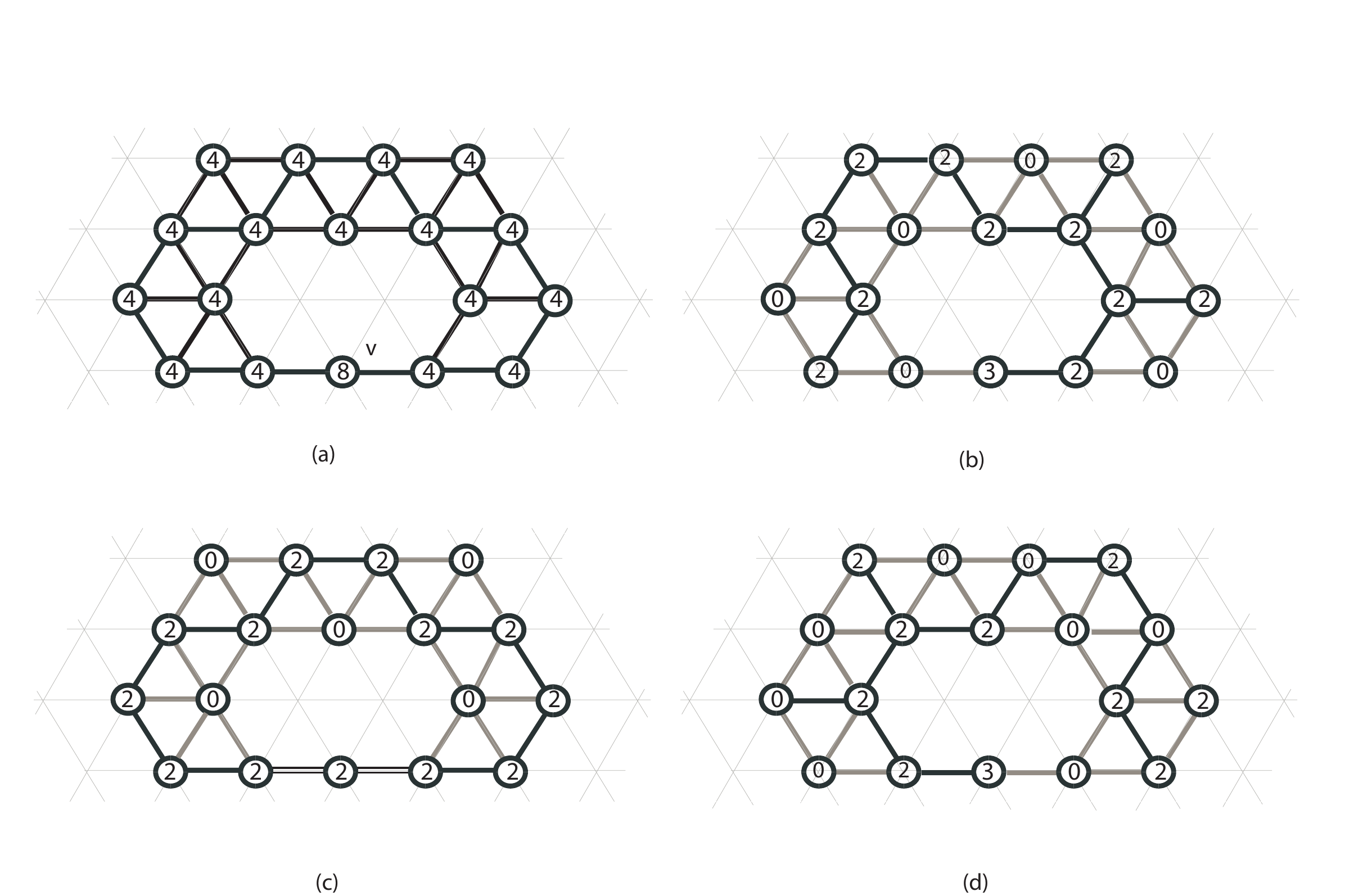}
 \end{center}
\vskip -5mm
 \caption{The example: a hexagonal graph (a) and  its  three triangle-free subgraphs (b), (c), and (d).} 
 \label{Example_decomposition}
\end{figure}

%%%%%%%%%%%%%%%%%%%%%%%%%%%%%%%%%%%%%%%%%%%%%%%%%%%%%%%%%%%%%%%%%%%%%%%%
%
\subsection{5/4 approximation algorithm  }

\begin{theorem} 
Let $G$ be a weighted hexagonal graph and assume all weights are even numbers.
There is a proper multicoloring of $G$ that uses at most   $15 \lfloor\frac{\omega(G)}{12}\rfloor+ 15 $     colors.
\end{theorem}

%\begin{proof}  
\proof
If    $\omega(G) < 12 $ then we can  
apply any 4/3-approximation algorithm to obtain a proper coloring of $G$ with at most 15 colors.
($\omega(G) \leq 11 $ gives $\frac{4 \omega(G) +1}{3} \leq 15$.)

Let $G  $ be a   weighted hexagonal graph with $\omega(G) \geq  12$. 
Define a   sequence  of  subgraphs of  $G_i$ as follows.
Start with $G_0 = G$.

In $G_i$, let $L_i \subseteq V_i$   be the set of vertices with low weight, $L_i = \{ v \in V_i  ~|~ d_i(v) < 4 \}$.
Define the subgraph  $\tilde G_i$ to be induced graph on $\tilde V_i$, where $\tilde V_i = V_i - L_i$ is the subset of $V_i$ of vertices with $d_i(v) \geq 4$. 
For a  vertex $v$ which is isolated in $\tilde G_i$  set $\tilde d_i(v) = 12$ and let  $\tilde d_i(v)$ equal 4 on vertices that lie on a triangle.  
For vertices that are not isolated and not on a triangle,  let 

$\tilde d_i(v) = 4$ if  $d_i(v) \leq  \frac{\omega(G_i)}{3}  $,

$\tilde d_i(v) = 6$ if  $  \frac{\omega(G_i)}{3}    <  d_i(v)  <   \frac{2\omega(G_i)}{3}   $,

$\tilde d_i(v) = 8$ if  $d_i(v) \geq  \frac{2\omega(G_i)}{3}  $.

Observe  that   $\tilde G_i$ can be colored with 15 colors. 
This
follows from   Lemma \ref{lemma12}
  because, by construction,  all weights of vertices in  $\tilde G_i$  are  either 4,6,8, or 12, and   $\omega(\tilde G_i) =12$.  

Define $G_{i+1}$  as follows.
For $v \in   \tilde V_i$   let     $d_{i+1}(v) =   d_i(v) - \tilde d_i(v)$  
(reduce the weights of vertices in $L_i \subseteq V_i$ and leave the  demands   of other vertices   unchanged). 

\noindent{\bf Claim.}  
 $\omega (G_{i+1}) \leq   \omega (G_i)  - 12  $ provided   $\omega (G_i) \geq 12$.

To see correctness of the  claim, assume  $ \omega(G_i ) \geq  12$.
The isolated vertices and triangles  of  $G_i$  obviosly get 12 colors, hence the weight of these cliques is reduced by 12.
Now consider a  maximal    clique  in  $G_i$   that  has exactly two  vertices $u$ and $v$, 
i.e. $    d_i(u) +  d_i(v) = \omega(G_i)$.
Distinguish three cases.
(1)    $\tilde d_i(u) = 4$.   Then    $  d_i(u)  \leq  \frac{\omega(G_i)}{3} $, and therefore    $d_i(v) \geq  \frac{2\omega(G_i)}{3}  $, 
and, by definition,     $\tilde d_i(v) =8$.  
(2)  $\tilde d_i(u) = 6$. 
 From   $  \frac{\omega(G_i)}{3}   <  d_i(v)  <   \frac{2\omega(G_i)}{3}   $ and  $d_i(u) +  d_i(v) =  \omega(G_i)$
we  get   $  \frac{\omega(G_i)}{3}      <  d_i(u) $ and $ d_i(u)  <   \frac{2\omega(G_i)}{3}   $, hence       $\tilde d_i(v) = 6$. 
(3)   Finally,  $\tilde d_i(u) = 8$  implies  $\tilde d_i(v) =4$.  

For  the largest cliques on two vertices   that are not maximal, we have    $    d_i(u) +  d_i(v) = \omega(G_i) -2$ because the demands are assumed to be even.
Observe that the    weights  of these cliques are   reduced by at least 10. 
This is true because it  is not possible to have   $   \tilde d_i(u)  =  \tilde d_i(v) = 4$  since $d_i(v) \leq  \frac{\omega(G_i)}{3}  $ and $d_i(u) \leq  \frac{\omega(G_i)}{3}  $
would imply   $     d_i(u)   +    d_i(v)    \leq  \frac{\omega(G_i)}{3}   $   
contradicting the fact that   $ \omega(G_i) -2 > \frac{\omega(G_i)}{3}   $ for $\omega (G_i)>   12$.

Clealy, the cliques on two vertices    for which     $    d_i(u) +  d_i(v) = \omega(G_i) -  4 =8 $   get  at least 4+4=8  colors.

Hence  $\omega (G_{i+1}) \leq   \omega (G_i)  - 12  $   as claimed.

Repeat  application of  the construction above to each   $G_i$ to obtain $G_{i+1}$,     until $\omega(G_i) < 12$.
From considerations above it follows that 
(1) the number of iterations is at most  $\omega(G)  \bdiv  12$.
(2) the total number of colors needed is at most    $15  (\omega(G) \bdiv  12 )$.
Hence for   $ k = \omega(G) \bdiv  12  =   \lfloor\frac{\omega(G)}{12}\rfloor  $ iterations we need   $15  k$  colors. 

Finally,  we have to color the last  graph $G_k$ with $\omega(G_k) < 12$  for which 15 colors suffice.
\KONEC

\begin{theorem}\label{theorem54}
Let $G$ be a weighted hexagonal graph.
There is a proper multicoloring of $G$ that uses at most   $15 \lfloor\frac{\omega(G)}{12}\rfloor+ 18 $     colors.
\end{theorem}

\proof
Reduce  the odd  weights of vertices of $G$ by one, and apply the previous theorem.
There may  be some vertices left that need one additional color.
Use three additional colors  (for example, the base colors) to fulfill the missing demand. 
\KONEC

\section{Conclusions and future work}
 
The partition of a weighted hexagonal graph  given in Subsection \ref{partition}
was used to provide a 15 coloring of a class of graphs with clique number 12.
This in turn enabled a 5/4  approximation for general weighted hexagonal graphs.

The present author believes that the mail idea, partitioning of hexagonal graph to triangle-free  subgraphs,  may be used to improve the bound 5/4 further.
There are 7/6 approximation  algorithms for triangle-free hexagonal  graphs, however they can not be 
applied to the partition of weights as proposed in this paper.  
It is an interesting question whether it is   possible to define  another,   more suitable  partition. 

More generaly, an interesting open question is whether it is possible to 
partition the weights of a hexagonal graph $G$  to   three  triangle-free subgraphs  $ G_1$, $G_2$, and $ G_3$ such that 
$\omega(G_1) + \omega(G_2) + \omega(G_3)  = \omega(G)$ ?
 Recall that the positive answer would imply that
given any  algorithm   that gives a $r$-approximation  for triangle-free hexagonal graphs  we would have a   $r$-approximation   for general hexagonal graphs !

\noindent
{\bf Acknowledgement.}
The author wishes to thank to Sundar   Vishwanathan  for spoting a gap  in the argument  in  the previous version of the manuscript \cite{preprint}.


\begin{thebibliography}{99}     \frenchspacing \small 

%\bibitem{a139} Chin, F.Y.L., Zhang, Y., Zhu H. \emph{A 1-local 13/9-competitive Algorithm for Multicoloring Hexagonal Graphs}, Algorithmica, vol. 54, pp 557-567 (2009)
\bibitem{Hale} Hale, W.K. \emph{Frequency assignment: theory and applications}, Proceedings of the IEEE, vol 68(12), pp 1497-1514 (1980) 
\bibitem{Havet} Havet, F. \emph{Channel assignment and multicoloring of the induced subgraphs of the triangular lattice}, Discrete Mathematics, vol. 233, pp 219-231 (2001) 
\bibitem {havetJZ} Havet, F., \v{Z}erovnik, J. \emph{Finding a Five
Bicolouring of a Triangle-free Subgraph of the Triangular Lattice}, Discrete
Mathematics vol. 244, pp 103-108 (2002)
%\bibitem{Janssen} Janssen, J., Krizanc, D., Narayanan, L., Shende, S. \emph{Distributed Online Frequency Assignment in Cellular Network}, Journal of Algorithms, vol. 36(2), pp 119-151 (2000) 

\bibitem{Klostermeyer1}
 Klostermeyer, W.,  Zhang, CQ.
(2 + $\eps$)-Coloring of planar graphs with large odd-girth. Journal of Graph Theory 33(2): 109-119 (2000)

\bibitem{Klostermeyer2} 
 Klostermeyer, W.,  Zhang, CQ.
n-Tuple Coloring of Planar Graphs with Large Odd Girth. Graphs and Combinatorics 18(1): 119-132 (2002)

\bibitem{pierwsza} McDiarmid, C., Reed, B. \emph{Channel assignment and weighted coloring}, Networks,  vol. 36(2), pp. 114-117 (2000) 
%\bibitem{channel} Narayanan, L. \emph{Channel assignment and graph multicoloring}, Handbook of wireless networks and mobile computing, pp 71-94, Wiley, New York, (2002) 
\bibitem{cykle} Narayanan, L., Shende, S.M. \emph{Static frequency assignment in cellular networks}, Algorithmica, vol. 29(3), pp 396-409 (2001) 
\bibitem{ost} Sau, I., \v{S}parl, P., \v{Z}erovnik, J. \emph{7/6-approximation Algorithm for Multicoloring Triangle-free Hexagonal Graphs}, Discrete Mathematics, vol. 312, pp 181-187 (2012)
\bibitem{a73lin} \v{S}parl, P., Witkowski, R. \v{Z}erovnik, J. \emph{A Linear Time Algorithm for $7-[3]$-coloring Triangle-free Hexagonal Graphs}, Information Processing Letters, vol. 112, pp. 567-571 (2012)
\bibitem{clanek3D}   \v{S}parl, P., Witkowski, R. \v{Z}erovnik, J., \emph{ Multicoloring of cannonball graphs},  Ars Mathematica Contemporanea,  vol. 10,  31 - 44   (2016) 
\bibitem{a54} \v{S}parl, P., \v{Z}erovnik, J. \emph{2-local 5/4-competitive algorithm for multicoloring triangle-free hexagonal graphs}, Information Processing Letters, vol. 90(5), pp 239-246 (2004) 
\bibitem{a43} \v{S}parl, P., \v{Z}erovnik, J. \emph{2-local 4/3-competitive Algorithm for Multicoloring Hexagonal Graphs}, Journal of Algorithms, vol. 55(1), pp 29-41 (2005) 
\bibitem{a76subclass} \v{S}parl, P., \v{Z}erovnik, J. \emph{2-local 7/6-competitive algorithm for multicoloring a sub-class of hexagonal graph},  International Journal of Computer Mathematics, vol. 87, pp 2003-2013 (2010)
\bibitem{tfree} Sudeep, K.S., Vishwanathan, S. \emph{A technique for multicoloring triangle-free hexagonal graphs}, Discrete Mathematics, vol. 300, pp. 256-259 (2005) 
%\bibitem{a1712} Witkowski, R. \emph{A 1-local 17/12-competitive Algorithm for Multicoloring Hexagonal Graphs}, Lecture Notes of Computer Science, vol. 5699/2009, pp 346-356 (2009)
\bibitem{a43moj} Witkowski, R. \emph{1-local 4/3-competitive Algorithm for Multicoloring a Subclass of Hexagonal Graphs},  Discrete Applied Mathematics vol. 164, pp 349-355 (2014) 
\bibitem{Utilitas} Witkowski, R. \v{Z}erovnik, J. \emph{Proof of McDiarmid-Reed conjecture for a subclass of hexagonal graphs}, Utilitas Mathematica, to appear.
%\bibitem{a65} \v{Z}erovnik, J. \emph{A distributed 6/5-competitive algorithm for multicoloring triangle-free hexagonal graphs}, International Journal of Pure and Applied Mathematics, vol. 23(2), pp 141-156 (2005) 
\bibitem{preprint} \v{Z}erovnik, J. \emph{Improved approximative multicoloring of hexagonal graphs}, 	arXiv:1606.01328v1 [math.CO]. 
\end{thebibliography}
\end{document}